\documentstyle[12pt]{article}
\begin{document}
\newtheorem{guess}{Proposition }[section]
\newtheorem{theorem}[guess]{Theorem}
\newtheorem{lemma}[guess]{Lemma}
\newtheorem{corollary}[guess]{Corollary}
\def \ja {\vrule height 3mm width 3mm}
\def\inftwo#1#2#3{\smash{\mathop{#1}\limits_{#2\atop{#3}}}}

\centerline{\Large \bf Moduli of bundles on the blown-up plane}
\vspace{.2 in}
\centerline{\Large Elizabeth  Gasparim
{\footnote{Partially supported
by a research grant from CNPQ}}}

\begin{abstract}{\small \em Let ${\cal M}_j$ denote the moduli space
of bundles on the blown-up plane which  
split as ${\cal O}(j) \oplus {\cal O}(-j)$ 
over the exceptional divisor.
We show that there is a topological embedding
$\Phi_j: {\cal M}_j \rightarrow {\cal M}_{j+1}.$ }
\end{abstract}

Let $E$ be a holomorphic bundle on the blown-up plane 
which restricts to the exceptional divisor as 
${\cal O}(j) \oplus {\cal O}(-j). $
Then there is a canonical form of transition matrix 
for  $E$ (see [1] Thm. 2.1). Namely,
for a good choice of coordinates
 there is a polynomial 
$p= \sum_{i = 1}^{2j-2}\sum_{l = i-j+1}^{j-1}p_{il}z^lu^i $ such that 
$E$ is defined  by 
the transition matrix 
$\left(\matrix{z^j & p \cr
                0 & z^{-j}}\right).$
(Note that in  particular it follows that 
$E$ is algebraic and is an extension of line bundles.)

Our aim here is to show that the map 
$\left(\matrix{z^j & p \cr
                0 & z^{-j}}\right)
\mapsto
\left(\matrix{z^{j+1} & z\,u^2p \cr
                0 & z^{-j-1}}\right)$
induces an inclusion of moduli spaces 
$\Phi_j: {\cal M}_j \rightarrow {\cal M}_{j+1}$
which is a homeomorphism over its image.

Fix an  integer $j.$ 
Then each polynomial  determines a rank two bundle over 
the blown-up plane.  
We write  $p \sim p'$ 
if the matrices $\left(\matrix{z^j & p \cr
                0 & z^{-j}}\right)$
and $\left(\matrix{z^j & p' \cr
                0 & z^{-j}}\right)$  
define isomorphic bundles.
Our 
polynomials have $N = (j-1)(2j-1)$ complex coefficients (see[1]).
Hence the space
${\bf C}^{N}/\sim$ 
(where two points are equivalent 
if the corresponding polynomials define
isomorphic bundles)
parametrizes  
isomorphism  classes of 
 bundles on 
the blown-up plane with splitting 
${\cal O}(j) \oplus {\cal O}(-j) $
over the exceptional divisor.
Denote the set of these isomorphism classes by $M_j$ 
and give ${\bf C}^{N}/\sim$
the quotient topology.
We have just shown that there
 is a natural identification of
${\bf C}^{N}/\sim$ with $M_j $ which therefore 
 inherits an induced 
topology ${\cal A}.$
We denote by ${\cal M}_j$ the topological space $(M_j, {\cal A}).$

\begin{guess}: The following map is an
embedding (i.e. homeomorphism over its image)
$$\begin{array}{rcl}\Phi_j: {\cal M}_j & \rightarrow &  {\cal M}_{j+1} \cr
 p & \mapsto &  z\,u^2p. \end{array} $$

\end{guess}

\noindent{\bf Proof:}
We first show that the map is well defined. 
Suppose 
$\left(\matrix{ z^j & p \cr 0 & z^{-j} }\right)$ and 
$\left(\matrix{ z^j & p' \cr 0 & z^{-j} }\right)$ 
represent isomorphic bundles. 
Then there are coordinate changes
 $ \left(\matrix{a  & b   \cr c  & d \cr}\right)$
holomorphic in $z, \, u$ and 
$ \left(\matrix{\alpha & \beta  \cr \gamma & \delta \cr}\right)$
  holomorphic in $z^{-1}, \, zu$
for which the following equality holds
(compare [1] pg. 587)
$$ \left(\matrix{\alpha & \beta  \cr \gamma & \delta \cr}\right)
 = \left(\matrix{z^j & p'  \cr 0 & z^{-j} \cr}\right)
 \left(\matrix{a  & b   \cr c  & d \cr}\right)  
\left(\matrix{z^{-j} & -p  \cr 0 & z^j \cr}\right).$$
Therefore these two bundles are isomorphic 
exactly when the system of equations 
$$ \left(\matrix{\alpha & \beta  \cr \gamma & \delta \cr}\right)
= \left(\matrix{a + z^{-j}p'c & 
z^{2j}b +z^j(p'd  -ap) - pp'c  \cr
z^{-2j}c & d - z^{-j}pc \cr}\right) \eqno (*)$$
can be solved by a matrix 
 $ \left(\matrix{a  & b   \cr c  & d \cr}\right)$
 holomorphic in $z, \, u$ which makes 
$ \left(\matrix{\alpha & \beta  \cr \gamma & \delta \cr}\right)$
 holomorphic in $z^{-1}, \, zu.$

On the other hand, the images of these two bundles are given by
transition matrices
$\left(\matrix{ z^{j+1} & z\,u^2p \cr 0 & z^{-j-1} }\right)$ and 
$\left(\matrix{ z^{j+1} & z\,u^2p' \cr 0 & z^{-j-1} }\right),$
which  represent isomorphic bundles iff there are
coordinate changes
  $ \left(\matrix{\bar a  & \bar b   \cr \bar c  & \bar d \cr}\right)$
holomorphic in $z, \, u$ and 
$ \left(\matrix{\bar\alpha & \bar\beta  \cr \bar\gamma & 
\bar\delta \cr}\right)$
  holomorphic in $z^{-1}, \, zu$ 
satisfying the equality
$$ \left(\matrix{\bar\alpha & \bar\beta  \cr \bar\gamma & 
\bar\delta \cr}\right)
 = \left(\matrix{z^{j+1} & z\,u^2p'  \cr 0 & z^{-j-1} \cr}\right)
 \left(\matrix{\bar a  & \bar b   \cr \bar c  & \bar d \cr}\right)  
\left(\matrix{z^{-j-1} & -z\,u^2p  \cr 0 & z^{j+1} \cr}\right).$$
That is, the images represent isomorphic bundles if the system
$$ \left(\matrix{\bar\alpha & \bar\beta  \cr \bar\gamma & 
\bar\delta \cr}\right)
= \left(\matrix{\bar a + z^{-j } u^2p'\bar c & 
z^{2j+2}\bar b +z^{j+2}u^2(p'\bar d  -\bar a p) - z^2u^4pp'\bar c  \cr
z^{-2j-2}\bar c & \bar d - z^{-j} u^2p\bar c \cr}\right) \eqno (**)$$
has a solution.

 Write $x = \sum x_i u^i$ for $x \in \{a,b,c,d,
 \bar a,\bar b,\bar c,\bar d\}$
and choose
$\bar a_i = a_{i+2},$ 
$\bar b_i = b_{i+2}u^2,$
$\bar c_i = c_{i+2}u^{-2},$
$\bar d_i = d_{i+2}.$  
Then if $ \left(\matrix{a  & b   \cr c  & d \cr}\right)$ 
solves  (*), one verifies  that 
 $ \left(\matrix{ \bar a  & \bar b   \cr \bar c  & \bar d \cr }\right)$
solves 
(**), which implies that the images  represent 
 isomorphic bundles and therefore $ \Phi_j $  is well defined.
To show that the map is injective just reverse the 
previous argument.
Continuity is obvious. 
Now we observe also that the image 
 $\Phi_j({\cal M}_j)$ is a saturated set in ${\cal M}_{j+1}$
 (meaning that if 
$ y \sim x$ and $x \in \Phi_j({\cal M}_j)$ then $y \in \Phi_j({\cal M}_j$)).
In fact, if $E \in 
 \Phi_j({\cal M}_j)$ 
then $E$ 
splits in the 2nd formal  neighborhood.
Now if $E' \sim E$ than $E'$ must also
split in the 2nd formal  neighborhood
therefore the polynomial corresponding to $E'$ is 
of the form $u^2p'$
and hence $ \Phi_j(z^{-1}p') $ gives $E'.$
Note also that 
 $\Phi_j({\cal M}_j)$ is a closed 
subset of ${\cal M}_{j+1},$ given by
the equations 
$p_{il} = 0$ for $i = 1,2$ and  $i-j+1 \leq l \leq j-1.$
Now the fact that  $\Phi_j$ is a homeomorphism over its  
image follows from the following  easy lemma.

\hfill\ja

\begin{lemma} 
Let $X \subset Y$ be a closed subset 
and $\sim$ an equivalence relation in $Y,$ such 
that $X$ is $\sim$ saturated. 
Then the map $I :X/{\sim} \rightarrow Y/{\sim}$
induced by the inclusion is a homeomorphism over the image.
\end{lemma}

\noindent{\bf Proof:} Denote by 
$\pi_X : X \rightarrow X/{\sim} $
and $\pi_Y: Y \rightarrow Y/{\sim}$ the projections.
Let $F$ be a closed subset of $X/{\sim}.$ Then 
$\pi_X^{-1}(F)$ is closed and saturated in $X$ and therefore 
$\pi_X^{-1} (F) $ is also closed and saturated in $Y.$ 
It follows that $\pi_Y(\pi_X^{-1}(F))$  is closed 
in $Y/{\sim}.$ \hfill\ja

\vspace{5 mm}

\noindent Acknowledgments:
This work was written during a visit to   
the International Centre for Theoretical Physics in Trieste.
I am thankful to Prof. M.S. Narasimhan for 
once more making it possible for me to visit 
the ICTP.

\noindent {\small Address:}

\noindent{\small Departamento de Matem\'atica, Universidade Federal de Pernambuco} 

\noindent {\small Cidade Universit\'aria, Recife, PE, BRASIL, 50670-901}

\noindent{\small gasparim@dmat.ufpe.br}

\end{document}